\documentclass[11pt,english]{revtex4-1}
\usepackage{mathptmx}
\usepackage[T1]{fontenc}
\usepackage[latin9]{inputenc}
\usepackage{babel}
\usepackage{amsmath}
\usepackage{amssymb}
\usepackage{graphicx}
\usepackage{esint}
\usepackage[unicode=true,pdfusetitle,
 bookmarks=true,bookmarksnumbered=false,bookmarksopen=false,
 breaklinks=false,pdfborder={0 0 1},backref=false,colorlinks=false]
 {hyperref}
\usepackage{breakurl}
\begin{document}

\title{Large Deviations Theory of Increasing Returns}

\author{Simone Franchini and Riccardo Balzan}

\affiliation{Sapienza Universit� di Roma, Piazza A. Moro 1, 00185 Roma, Italy}
\begin{abstract}
An influential theory of increasing returns has been proposed by the
economist W. B. Arthur in the '80s to explain the lock-in phenomenon
between two competing commercial products. In the most simplified
situation there are two competing products that gain customers according
to a majority mechanism: each new customer arrives and asks which
product they bought to a certain odd number of previous customers,
and then buy the most shared product within this sample. It is known
that one of these two companies reaches monopoly almost surely in
the limit of infinite customers. Here we consider a generalization
{[}G. Dosi, Y. Ermoliev, Y. Kaniovsky, J. Math. Econom. 23, 1\textendash 19
(1994){]} where the new customer follows the indication of the sample
with some probability, and buy the other product otherwise. Other
than economy, this model can be reduced to the urn of Hill, Lane and
Sudderth, and includes several models of physical interest as special
cases, like the Elephant Random Walk, the Friedman's urn and other
generalized urn models. We provide a large deviation analysis of this
model at the sample-path level, and give a formula that allows to
find the most likely trajectories followed by the market share variable.
Interestingly, in the parameter range where the lock-in phase is expected,
we observe a whole region of convergence where the entropy cost is
sub-linear. We also find a non-linear differential equation for the
cumulant generating function of the market share variable, that can
be studied with a suitable perturbations theory.
\end{abstract}
\maketitle
\tableofcontents{}

\part{Main results}

\section{Introduction\label{sec:Introduction}}

It is known that certain economic markets - especially the technological
ones - show \textit{increasing returns} \cite{Arthur nature,Arthur book,Arthur},
a positive feedback phenomenon where if a company gains some initial
advantage (even small) is more likely to get even more in the future,
eventually dominating the market share in the long run - this phase
is also called \textit{lock-in} into a monopolistic state. To understand
the origin of this effect, a simplified market model has been introduced
in the '80s by the economist W. B. Arthur \cite{UM Arthur} in the
framework of its Increasing Returns theory (IRT) \cite{Arthur book,Arthur}.
Let consider two competing companies that launch a new kind of product
roughly at the same time (as practical example we could think about
two smart phones in the early 2000s). Suppose that these products
are roughly equivalent, such that there is no practical reason for
choosing one over the other, we can imagine that a buyer will base
his decision in part on personal opinions (personal tastes, ideologies,
advertising, etc.) and in part on those of other people that already
purchased one of the products. Then, let us consider a simplified
situation in which the new customers are imperfectly informed about
the products, so that they will make their choices by looking at the
number of adopters who are already using it \cite{Arthur,UM Arthur,Dosi Ermoliev,Dosi last}.
An alternative hypothesis that gives the same effect is to consider
positive (or negative) externalities in adoption \cite{Dosi Ermoliev,Dosi last}.
In both cases, we consider the additional rule that any new adopter
will choose the technology used by the majority of the sample only
with a certain probability, and the other technology otherwise \cite{Dosi Ermoliev,Dosi last}. 

This scenario has been considered by G. Dosi, Y. Ermoliev and Y. Kaniovsky
(DEK, 1994) \cite{Dosi Ermoliev}, the proposed model is as follows:
consider a binary vector that represents the individual choices of
the customers, 
\begin{equation}
X_{N}:=\{X_{N,1},\,X_{N,2},\,...\,,X_{N,N}\},
\end{equation}
with $X_{N,n}\in\{0,1\}$ and $N$ potential size of the market. This
vector represents the full history of the market evolution, from the
first sell to full saturation, when the maximal number of customers
is reached. The variable $X_{N,n}$ represents the choice of the $n-$th
customer, we arbitrarily associate the value one to the first product
and zero to the second. The total number of customers of the first
product will therefore be
\begin{equation}
\Gamma_{N,n}:=\sum_{m\leq n}X_{N,m}
\end{equation}
the market share of the first product up to the $n-$th customer is
represented by the variable
\begin{equation}
x_{N,n}:=\frac{\Gamma_{N,n}}{n}=\frac{1}{n}\sum_{m\leq n}X_{N,m}.
\end{equation}
Then, the choice of the next customer $X_{N,n+1}$ is determined by
the following rule: first, sample $k$ previous customers, where $k$
is an odd integer (this to avoid inconclusive outputs from the poll).
Then, if the sample is found to have more customers that bought the
first product, the variable $X_{N,n+1}$ will be equal to one with
a probability $p$, and will be zero otherwise. On the other hand,
if more customers owning the second product are found in the sample,
$X_{N,n+1}$ will be zero with probability $p$, and one otherwise.
Notice that the new customers follow the majority of the polled sample
with probability $p$, that in some sense quantifies the trust of
the newcomers in the behavior of their predecessors: hereafter we
will call $p$ \textit{trust parameter}, although it may also reflect
more practical constraints, such as a requirement for compatibility
with the technology adopted by the polled customers. For $p=1$ the
DEK model describes a market where the customers always buy the product
owned by the majority of the sample, i.e., the original version introduced
by Arthur et al. (1983) and Arthur (1989) \cite{UM Arthur,Arthur}:
some sample trajectories of this process for $p=1$ and $k=3$ are
in Figure 2a of G. Dosi et al. (2017) \cite{Dosi last}. Concerning
the initial conditions, we will distinguish of two kinds: we introduce
$\tau\in\left[0,1\right]$ the fraction of customers that made their
choices already (\textit{market saturation} parameter): in this paper
we will consider an \textit{early start} in the market at some fixed
number of customers $M<\infty$, also called \textit{virgin market}
condition, that in the limit of infinite customers is equivalent to
a debut in the market approximately at $\tau=0$ (and does not affect
the LDT theory for $N\rightarrow\infty$) and a \textit{late start}
$M=\tau N$ (a product that enters in the market when the saturation
is already macroscopic), that strongly influences the distribution
of the final share also at the LDT level.

\section{Relation with HLS urns}

In this paper we develop a Large Deviations theory (LDT) for the DEK
model for any $p$ and $k$ by adapting results from the Hill Lane
and Sudderth (HLS) urn model \cite{HLS1,HLS2}, a very general model
for which a mathematically rigorous LDT has been recently developed
\cite{Franchini}, and that includes the DEK model as special case.
An HLS urn process \cite{Pemantle review,Mahmoud,HLS1,HLS2,Pemantle Touch,Franchini}
is a two color urn process controlled by a functional parameter $\pi\left(x\right)$
that we call \textit{urn function} (actually \textit{adoption function}
in Ref. \cite{Arthur}), where the new step $X_{N,n+1}$ is one with
probability $\pi\left(x_{N,n}\right)$ and zero otherwise. The relation
between IRT and HLS urns is well known since the very beginning, in
fact, this model has been introduced independently by HLS (1980) and
then also by Arthur et al. (1983) within just three years. The urn
function that describes the DEK model can be determined as follows:
start with $k=1$, the probability of extracting an owner of the first
product is $x$, then their total number will increase with probability
\begin{equation}
\pi_{1}\left(x\right):=p\,x+\left(1-p\right)\left(1-x\right)=\left(1-p\right)+\left(2p-1\right)x,
\end{equation}
that is a linear urn function. In case $k=3$: the probability of
increasing the owners of the first product is that of extracting two
positive and one negative, plus that of extracting three positive,
that is \cite{Dosi Ermoliev} 
\begin{equation}
P_{3}\left(x\right):=x^{3}+3x^{2}\left(1-x\right)=3x^{2}-2x^{3},
\end{equation}
then, the corresponding urn function is \cite{Dosi Ermoliev}
\begin{equation}
\pi_{3}\left(x\right):=p\,\left(3x^{2}-2x^{3}\right)+\left(1-p\right)\left(1-\left(3x^{2}-2x^{3}\right)\right)=\left(1-p\right)+\left(2p-1\right)\left(3x^{2}-2x^{3}\right)\label{eq:pik-1}
\end{equation}
and cannot be reduced to the linear case $k=1$. In general, the probability
of finding a positive majority when extracting an odd number $k$
of steps is \cite{Dosi Ermoliev} 
\begin{equation}
P_{k}\left(x\right):=\sum_{h>k/2}\frac{k!}{h!\left(k-h\right)!}\,x^{h}\left(1-x\right)^{k-h}
\end{equation}
where the $h$ sum runs from $\left(k+1\right)/2$ to $k$. Follows
that the urn function that describes a DEK model with $k>2$ extractions
per step is \cite{Dosi Ermoliev} 
\begin{equation}
\pi_{k}\left(x\right):=p\,P_{k}\left(x\right)+\left(1-p\right)\left(1-P_{k}\left(x\right)\right)=\left(1-p\right)+\left(2p-1\right)P_{k}\left(x\right),\label{eq:pik}
\end{equation}
this is a $k-$th degree polynomial, and is therefore non-linear for
all non-trivial values of the trust parameter $p$. 

\begin{figure}
\includegraphics[scale=0.29]{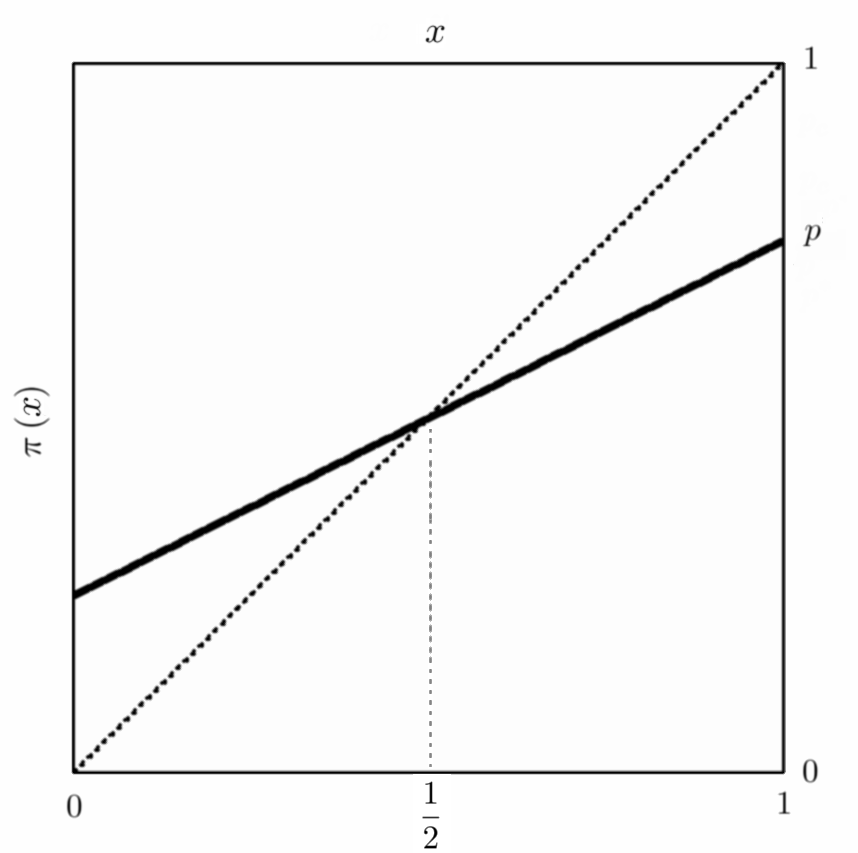}

\caption{\label{fig1}Example of linear urn function $\pi_{1}\left(x\right)$
for the DEK with $k=1$, the memory parameter is $p=5/8$. The urn
function always down-crosses the diagonal at $x_{0}=1/2$, that is
the only convergence point.}
\end{figure}

In case of a virgin market start, the convergence properties of the
HLS urns with any continuous urn functions have been studied in \cite{HLS1,HLS2,UM Arthur,Arthur,Dosi Ermoliev,Pemantle Touch,Franchini},
finding that the points of convergence of $x_{N,N}$ always belong
to the set of solutions of 
\begin{equation}
\pi\left(x\right)=x,\label{eq:urnqe}
\end{equation}
and that these solutions are stable only if the derivative of the
urn function in those points is smaller than one, i.e., if the $\pi\left(x\right)$
crosses $x$ from top to bottom (down-crossing). For the DEK model
with $k=1$, the urn function $\pi_{1}\left(x\right)$ crosses $x$
at $1/2$ for any value of $p<1$, and therefore $1/2$ is the only
possible point of convergence for the associated share $x_{N,N}$,
see Figure \ref{fig1}. This imply that $x_{N,N}$ converges to $1/2$
almost surely 
\begin{equation}
\lim_{N\rightarrow\infty}x_{N,N}=1/2,\ a.s.
\end{equation}
for all values of $p<1$ and of the initial condition $x_{N,M}$ a
phase diagram for the DEK $k=1$ is shown in Figure \ref{fig:3}.
This model does not show the lock-in phenomenon, although there is
still a value of $p$ where the dynamics is expected to slow down
(see Section \ref{sec:Cumulant-Generating-Function}). 

\begin{figure}
\includegraphics[scale=0.29]{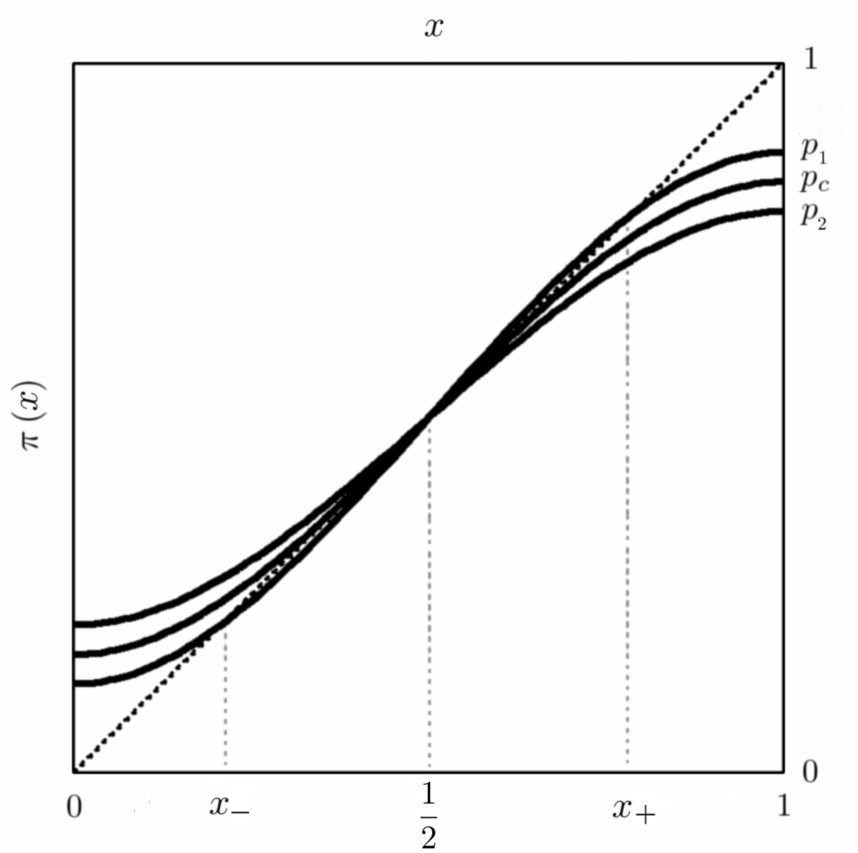}\caption{\label{fig:2}Three examples of the urn function $\pi_{3}\left(x\right)$
for a generalized DEK with $k=3$. The figure shows the urn functions
for three non-trivial memory parameters, $p_{2}=19/24$, $p_{c}=5/6$
and $p_{1}=21/24$. Below $p_{c}$ the urn function down-crosses the
line $x$ at $x_{0}=1/2$, that is the only convergence point. For
$p>p_{c}$ the point $x_{0}$ becomes an up-crossing (unstable equilibrium),
and the urn function crosses the diagonal $x$ also in $x_{-}$ and
$x_{+}$, that are both down-crossings and are the new stable attractors
for the process $x_{N}$. }
\end{figure}

\begin{figure}
\includegraphics[scale=0.29]{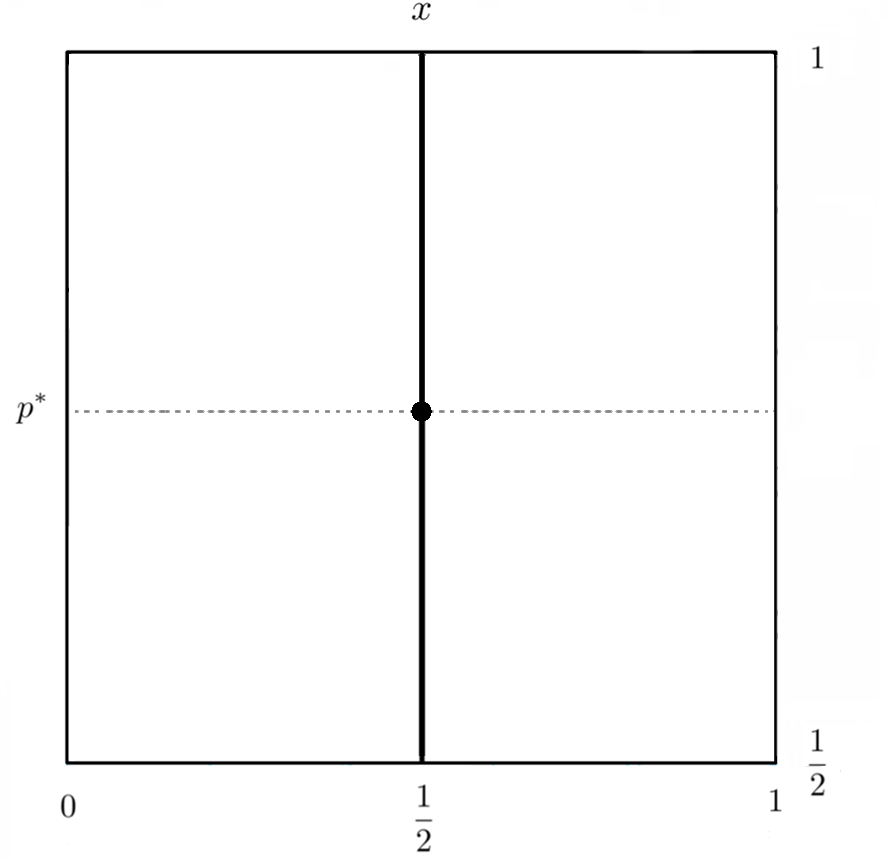}

\caption{\label{fig:3}Phase diagram $x$ vs $p$ for the entropy density $\phi\left(x\right)$
of the DEK $k=1$, the diagram is shown for $p>1/2$. For all $p<1$
the point $x_{0}=1/2$ is the only point of convergence for the density
of black balls, although there still is a critical $p^{*}$ where
the derivative of $\pi\left(x\right)$ in $x_{0}$ crosses the value
$1/2$, and the convergence of $x_{N}$ is slowed according to the
Pemantle mechanism \cite{Jack Harris,Pemantle Touch,Franchini}, see
Section \ref{sec:Perturbations-theory}. The line $x_{0}=0$ is always
a stable attractor for $x_{N}$, and the entropy is convex and strictly
negative in the whole region, except at the critical line $x=0$,
where is zero. According to Eq. (\ref{eq:urnqe}), there is a critical
value at $p^{*}=3/4$ at which the derivative of the urn function
gets above $1/2$: for $p>p^{*}$ there is a shape change in $\phi\left(x\right)$
in the neighborhood of $x=1/2$.}
\end{figure}

In the DEK with $k>2$ one can see the appearance of the lock-in phase
above some critical $p_{c}$. For $k=3$ the Eq. (\ref{eq:urnqe})
is a third degree equation, and can be solved with the well known
formula. In general, we find tree solutions \cite{Dosi Ermoliev}:
$x_{0}=1/2$ and 
\begin{equation}
2x_{\pm}=1\pm\sqrt{\frac{6p-5}{2p-1}},\label{eq:PARAK3}
\end{equation}
the quantity inside the square root is positive for $p\leq1/2$ and
$p\geq5/6$, but notice that $x_{\pm}\in\left[0,1\right]$ only if
$p\in\left[1/2,1\right]$, then, for $p$ below the critical value
$p_{c}=5/6$ there is again a unique stable solution at $1/2$ that
crosses $x$ from top to bottom (down crossing), see Figure \ref{fig:2}.
Above $p_{c}$ the function $\pi$ still crosses $x$ at the point
$1/2$, but it now does from bottom to top, i.e. it is an up-crossing
and is therefore not stable. Notice that for $p>p_{c}$ two new solutions
$x_{+}$ and $x_{-}$ appear, those are both down-crossings, and can
be stable attractors for $x_{N,N}$. Therefore, for $p$ above $p_{c}$
there are two attractors separated by an unstable equilibrium point
at $x_{0}=1/2$ \cite{Dosi Ermoliev}. Notice that in the limit of
infinite $k$ the probability of finding a majority of first product
owners within the sample converges to 
\begin{equation}
P_{\infty}\left(x\right):=\theta\left(1-2x\right)
\end{equation}
ie, the urn function $\pi_{k}$ converges to a step function
\begin{equation}
\pi_{\infty}\left(x\right):=\left(1-p\right)+\left(2p-1\right)\theta\left(1-2x\right)
\end{equation}
that still crosses the diagonal at the point $x_{0}=1/2$ (from top
to bottom) for $p<p_{c}=1/2$, and at $x_{-}=p$, $x_{+}=1-p$ if
the trust parameter is above $p_{c}$. Then also in the infinite $k$
limit there is a $p_{c}$ above which we find the same region of the
phase diagram that is observed for $k=3$. In fact, the phase diagram
shows the same structure for all $k>2$, apart from different $p_{c}$
and $x_{\pm}$. For this reason, we will concentrate our analysis
to the cases $k=1$ and $k=3$. 


\begin{figure}
\includegraphics[scale=0.29]{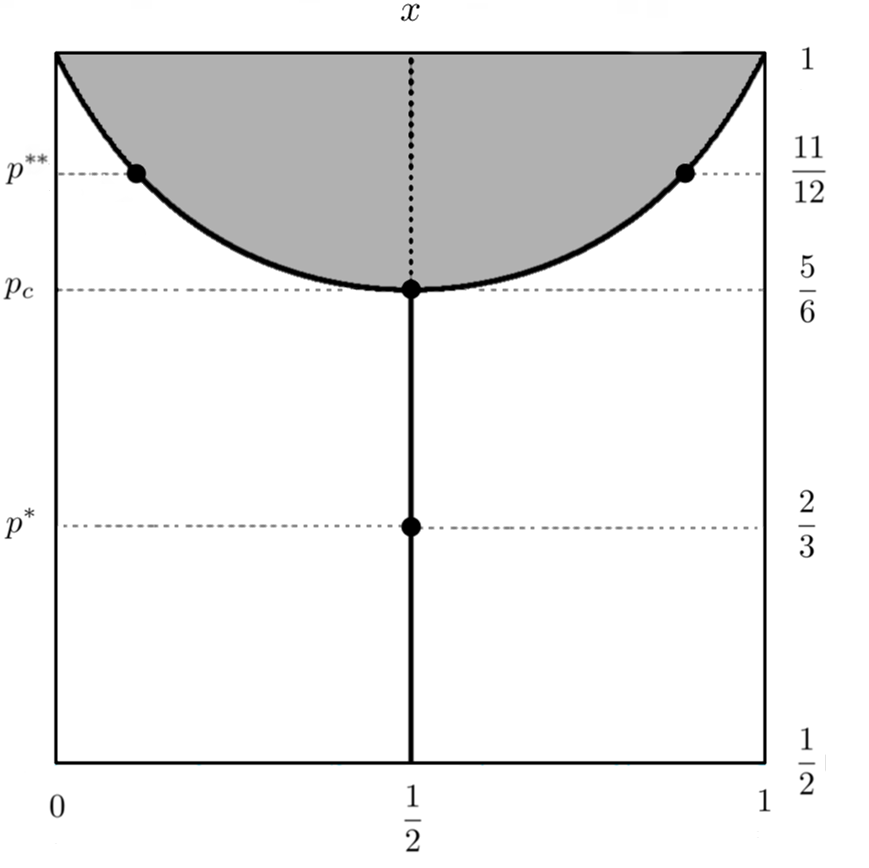}~\caption{\label{fig:4}Phase diagram $x$ vs $p$ for the entropy density $\phi\left(x\right)$
of the generalized DEK $k=3$. Above the critical value $p_{c}=5/6$
the point $x_{0}=0$ becomes an unstable equilibrium, and two new
symmetric attractors arise according to Eq. (\ref{eq:PARAK3}). In
the white colored region we still find a convex and negative $\phi\left(x\right)$,
except on the critical line, but notice that a new region appeared
above $p_{c}=5/6$, highlighted in darker shade, where $\phi\left(x\right)=0$,
i.e. the entropy is sub linear in $N$. The shape of $\phi\left(x\right)$
near the critical line is similar to the case $k=1$ for $p<p_{c}$,
except that here the point $p^{*}=2/3$ at which the derivative of
the urn function rise above $1/2$. Below $p_{c}$ the derivative
of $\pi$ in $x_{0}$ is increasing in $p$, and $p^{*}$is the value
at which crosses the value $1/2$ (from below). On the contrary, above
$p_{c}$ the derivative of $\pi$ in $x_{\pm}$ decreases in $p$,
and crosses the value $1/2$ (from above) at $p^{**}=11/12$: when
the derivative of the associated urn function in $y_{\pm}$ goes back
below $1/2$ the shape of $\phi\left(x\right)$ changes in the right
(left) neighborhood of $x_{-}$($x_{+})$. See also Figure \ref{fig:2}.}
\end{figure}

Summarizing, under virgin market condition the limit value of $x_{N,N}$
for $p>p_{c}$ converges to the points $x_{-}$ and $x_{+}$ almost
surely for any initial condition $x_{N,M}$ with $M<\infty$ (the
phases for $k=3$ are shown in the Figure \ref{fig:4}) but since
the urn functions that we are considering never touches zero or one
at any $x\in\left(0,1\right)$, for any initial condition $x_{N,M}$
that is fixed at $M\ll N$ there is a strictly positive probability
to reach the nearby of any other $x$ by gaining a finite number of
customers at the beginning of the process, then in the limit $N\rightarrow\infty$
both points $x_{\pm}$ carry some non-zero probability mass for any
early start. Anyway, it can be shown that the probability mass of
that point farther from $x_{N,M}$ will be exponentially suppressed
as $M$ grows. Fixing the initial condition at some $M=o\left(N\right)$
but still divergent in $N$ will suppress one of the two possibilities,
and concentrate the probability mass in the attraction point $x_{\pm}$
that is closest to the initial share $x_{N,M}$. Concerning the case
of late market entry at some $M=\tau_{0}N$, we discuss it in Section
\ref{sec:Main-results}, after introducing the optimal trajectories.

\section{Relation with other models}

We remark that, apart from economic models, the theory of the HLS
urns allows to put IRT in relation with many others interesting situations
that can be embedded (or approximated) by this very general urn model:
there is a number of computer science problems on preferential attachment,
network growth etc. \cite{Pemantle review,Mahmoud} that can be studied
in this framework. Here we list three that, in our opinion, are of
special physical interest. Transfer of knowledge between these field
would be certainly fruitful, and should be encouraged. 

For example, the case $k=1$ of the DEK is fully equivalent to another
well known stochastic model, the Elephant Random Walk (ERW), a simple
random walk where each new step is determined by selecting one of
the previous, then going in the same direction with probability $p$.
This model appears to have been re-descovered independently by G.
Sch�tz and S. Trimper in 2004 \cite{ERW shcutz trimper}, ten years
after the introduction of the DEK model, and has received much attention
since then as a paradigmatic example of processes with long range
memory. An important advancement in the understanding of this model
was made in 2016, when E. Baur and J. Bertoin observed \cite{ERW UM Baur Berton}
that the ERW could be mapped exactly into a two color urn of the Friedman's
type (that is in fact equivalent to a linear HLS urn \cite{Pemantle review,Mahmoud,Franchini})
where at each time one ball is drawn from the urn, and then replaced
together with a fixed numbers of new balls whose color depend on which
was drawn. This finding allowed many quantities of interest to be
studied from known results on these types of models, \cite{ERW UM Baur Berton,Jack Harris}
however, this analogy cannot be extended to the $k>2$ case. 

Also, Jack \cite{Jack LD} has identified the urn function describing
an interesting irreversible growth model introduced by Klymko, Garrahan
and Whitelam \cite{KGW,KGGW}, and also this model exhibits a lock-in
phase with a sub-linear entropy region, that is similar to the $k>2$
case of the DEK model \cite{Jack LD}. In this perspective, it would
be quite interesting to investigate also the universal HLS scaling
for symmetric urn functions recently proposed by Nakayama et al. (2021)
\cite{Kazuaki}. Notice that in Jack (2019) \cite{Jack LD} a non-rigorous
but powerful LDT is presented for a large class of models, whose predictive
power should be comparable to the rigorous LD techniques used in Franchini
(2017), see also the interesting review Jack 2020 \cite{Jack LD-1}.

Finally, the HLS framework allows to relate the DEK model with the
very classic Random Walk Range problem \cite{Huges,Franchini Range,Franchini Range Urns,Franchini Range Line},
that studies the number of different sites visited by a random walk
on the lattice $\mathbb{Z}^{d}$. This problem is important to polymer
physics as it exhibits a coil-globule transition at some critical
range density, and is in the same universality class of the Self-Avoiding
Walk above that value \cite{Franchini Range}. In Ref. \cite{Franchini Range Urns}
is shown that the Range problem can be exactly embedded in the HLS
model for some non-linear urn function at any $d$. For $d=2,3$ a
strongly non-linear urn function is observed (by numerical analysis),
but for $d\geq4$ the urn function gets surprisingly close to some
linear function in the self-avoiding walk-like region of large range
values, that would then be related to the DEK model with $k=1$. This
model also shows a sub-linear entropy region below some critical range,
as can be deduced also from a very detailed analysis of the ``moderate
deviations'' of the Wiener Sausage in the collapsed phase by M. van
den Berg, E. Bolthausen, F. Den Hollander (2001) \cite{van den Berg}.
Interestingly, they find cases of non-homogeneous optimal trajectories
with sub-linear entropy cost: we conjecture that in this collapsed
region the range undergoes a mechanism similar to that observed in
the lock-in phase (actually, the non-homogeneous zero-cost trajectories
that are described in Corollary 7 of Ref. \cite{Franchini}). 


\section{Zero-cost trajectories\label{sec:Zero-cost-trajectories}}

We perform a Large Deviations analysis for the HLS model at the sample-path
level, and adapt it to find the most probable trajectories taken by
the DEK model. Let $\tau\in\left[0,1\right]$ be the level of market
saturation (or the fraction of customers that already made their choice):
the \textit{optimal trajectories}, that we indicate with the symbol
$u$, are the scaling limit for $n/N\rightarrow\tau$ of the most
likely trajectories followed by the share variable $x_{N,n}$ of the
first product to reach some given final share $x$. These can be obtained
by solving the variational problem that is presentend in Section \ref{sec:Large-deviations}
(see also Theorem 4 of Ref. \cite{Franchini} for a full mathematical
derivation). Most interesting, we will show that for any initial condition
with positive saturation $\tau_{0}>0$ the scaling limit of the trajectory
taken by the share variable 
\begin{equation}
\lim_{N\rightarrow\infty}x_{N,\left\lfloor \tau N\right\rfloor }=:u\left(\tau\right)
\end{equation}
is non-degenerate for any starting share $u_{0}\in\left[0,1\right]$,
and can be found by inverting the following integral: 
\begin{equation}
\tau\left(u\right)=\tau_{0}\,\exp\int_{u_{0}}^{u}\frac{d\alpha}{\pi\left(\alpha\right)-\alpha}.\label{eq:tau-1-2-1-1-1}
\end{equation}
A crucial quantity of our analysis will be the scaling limit of the
entropy (logarithm of the probability) of $x_{N,N}$ converging to
some given $x$. Define the asymptotic limit of the entropy per customer
(hereafter we will call it \textit{entropy density}) 
\begin{equation}
\phi\left(x\right):=-\lim_{N\rightarrow\infty}\frac{1}{N}\log P\left(x_{N,N}=\left\lfloor xN\right\rfloor /N\right),
\end{equation}
informally, this is the scaling limit of the entropy respect to the
total number of customers, i.e. for a large number of customers the
probability of reaching a share $x$ is proportional to
\begin{equation}
P\left(x_{N,N}=\left\lfloor xN\right\rfloor /N\right)\sim\exp\left(-N\phi\left(x\right)\right).
\end{equation}

In the Section \ref{sec:Large-deviations}  we will show that the
shape of the limit entropy density $\phi$ can be linked to the trajectories
taken by the number of customers of the first product $\Gamma_{N,n}$
to reach its final value $\Gamma_{N,N}=\left\lfloor xN\right\rfloor $.
These trajectories, that we indicate with the symbol $\varphi$, are
the scaling limit $n/N\rightarrow\tau$ for the number of customers
of the first product, rescaled with the total number of customers
$N$, 
\begin{equation}
\lim_{N\rightarrow\infty}\Gamma_{N,\left\lfloor \tau N\right\rfloor }/N=:\varphi\left(\tau\right)
\end{equation}
this is related to the scaling limit of the share by the formula
\begin{equation}
u\left(\tau\right)=\varphi\left(\tau\right)/\tau.\label{eq:trasf}
\end{equation}
In Theorem 4 of Ref. \cite{Franchini} (see Section \ref{sec:Large-deviations}
of the present paper for an informal derivation) it is shown that
the limit entropy density of any HLS urn model with $\alpha-$H�lder
urn function $\pi$ is obtained trough the following LD principle:
let $C_{1}\left(\left[0,1\right]\right)$ be the set of absolutely
continuous function on $\left[0,1\right]$ (essentially, such that
the derivative exists almost everywhere) and let $Q\subset C_{1}\left(\left[0,1\right]\right)$
the subset of those functions with initial value zero, and such that
their derivative is positive but smaller than one ($1-$Lipschitz
function), 
\begin{equation}
Q:=\{\varphi\in C_{1}\left(\left[0,1\right]\right):\,\partial_{\tau}\varphi\left(\tau\right)\in\left[0,1\right],\,\varphi\left(0\right)=0\},\label{eq:dfdfd-2}
\end{equation}
also, let $Q\left(x\right)$ be the subset with final value $x$,
\begin{equation}
Q\left(x\right):=\{\varphi\in Q:\,\varphi\left(1\right)=x\}.\label{eq:dfdfd}
\end{equation}
Also, define the auxiliary function
\begin{equation}
L\left(\alpha,\beta\right):=\alpha\log\left(\beta/\alpha\right)+\left(1-\alpha\right)\log\left(\left(1-\beta\right)/\left(1-\alpha\right)\right),\label{eq:sfsf}
\end{equation}
then, the entropy density $\phi\left(x\right)$ can be computed by
solving the following variational problem: 
\begin{equation}
\phi\left(x\right)=\inf_{\varphi\in Q\left(x\right)}I\left(\varphi\right),\label{eq:THEVARIATION}
\end{equation}
with rate function defined as follows: 
\begin{equation}
I\left(\varphi\right):=-\int_{0}^{1}d\tau\,L\left(\partial_{\tau}\varphi\left(\tau\right),\pi\left(\varphi\left(\tau\right)/\tau\right)\right).\label{ratefunction}
\end{equation}

From this general result we where able to deduce a method to identify
those trajectories followed by the process $x_{N,n}$ that have a
sub-linear entropy cost, i.e., such that the entropy cost of following
the trajectory $u\left(\tau\right)$ is of the order $o\left(N\right)$:
it implies that the probability of following such a trajectory decays
sub-exponentially in the number of customers (actually as a power
law, see in Section \ref{sec:Scaling-of-the}), and not exponentially
fast as for those with a cost linear in $N$. Hereafter we will improperly
call these the \textit{zero-cost trajectories}, although their absolute
entropy cost is not exactly zero. It is shown that these zero-cost
trajectories can be deduced from the variational problem in Eq. (\ref{eq:THEVARIATION}),
with the additional constraint that the Lagrangian of the Eq. (\ref{ratefunction})
is exactly zero. In Corollary 6 of the Ref. \cite{Franchini} explicit
formulas are derived for those optimal trajectories ending in the
region where $\phi\left(x\right)=0$. Since for the HLS model the
function $L$ is a negative concave function, the condition $I\left(\varphi\right)=0$
implies that the trajectory $\varphi$ satisfies the equation 
\begin{equation}
L(\partial_{\tau}\varphi\left(\tau\right),\pi\left(\varphi\left(\tau\right)/\tau\right))=0,
\end{equation}
if this condition can be explicited in the variable $\partial_{\tau}\varphi\left(\tau\right)$,
then it provides the differential equation for the zero-cost trajectories.
Remarkably, since $L(\alpha,\beta)=0$ if and only if $\alpha=\beta$,
then the condition before reduces to the autonomous equation
\begin{equation}
\partial_{\tau}\varphi\left(\tau\right)=\pi\left(\varphi\left(\tau\right)/\tau\right),
\end{equation}
with final condition $\varphi\left(1\right)=x$. Applying the substitution
in Eq. (\ref{eq:trasf}) we obtain the equation for the scaling of
the share,
\begin{equation}
\frac{\partial_{\tau}u\left(\tau\right)}{\pi\left(u\left(\tau\right)\right)-u\left(\tau\right)}=\frac{1}{\tau}
\end{equation}
with final condition $u\left(1\right)=x$. This equation can be integrated
exactly, in the end one finds that the trajectories $u\left(\tau\right)$
can be computed in implicit form: the result is a simple formula,
\begin{equation}
\Pi\left(u\right)-\Pi\left(x\right)=\log\tau\left(u\right),\label{eq:tau-1}
\end{equation}
where $\Pi$ is the primitive function (indefinite integral) of the
reciprocal of $\pi\left(\alpha\right)-\alpha$, that is 
\begin{equation}
\Pi\left(\alpha\right):=\int\frac{d\alpha}{\pi\left(\alpha\right)-\alpha}.
\end{equation}
We can formally invert the formula of $\tau\left(u\right)$ before,
and write the equation for the zero-cost trajectories as follows:
let $\Pi^{-1}$ the inverse function of $\Pi$, then
\begin{equation}
u\left(\tau\right)=\Pi^{-1}\left(\Pi\left(x\right)+\log\tau\right).\label{eq:tau-1-2-1}
\end{equation}
The first important remark about this formula is that it allows to
extend the convergence theory of HLS urns also in case of a late start
in the market: let $\tau_{0}$ be the level of market saturation,
and let $u_{0}$ be the initial share (for a firm entering in the
market at $\tau_{0}$ we would have $u_{0}=0$), then, it can be shown
by inverting Eq. (\ref{eq:tau-1-2-1}) that for any positive $\tau_{0}$
there is a unique point 
\begin{equation}
\lim_{N\rightarrow\infty}x_{N,N}=x\left(u_{0},\tau_{0}\right)\ \ a.s.
\end{equation}
where the final share $x_{N,N}$ converges almost surely, 
\begin{equation}
x\left(u_{0},\tau_{0}\right):=\Pi^{-1}\left(\Pi\left(u_{0}\right)-\log\left(\tau_{0}\right)\right),\label{eq:tau-1-2-2}
\end{equation}
it can be shown that the convergence points found before for the virgin
market case are recovered in the limit $\tau_{0}\rightarrow0$.

In general, these results about the optimal trajectories could be
useful to confront with (and then eventually fit) trajectories followed
by real datasets in those cases where both the time series of the
share and the saturation are known for the considered market \cite{Dosi last}.
In this respect, it is important to realize that the market saturation
is not a time variable: the process describes the competition between
firms, but does not need to specify the underlying market grow. Let
$n\left(t\right)$ be the total number of customers up to time $t$,
growing according to some law in such way that the limit market size
is finite and equal to $N$. Then, the share of the first product
up to time $t$ would be $x\left(t\right)=x_{N,n\left(t\right)}$,
that can be confronted with the predicted scaling limit 
\begin{equation}
\lim_{N\rightarrow\infty}x\left(t\left(\tau\right)\right)=u\left(\tau\right).
\end{equation}
by plotting $x\left(t\right)$ in function of the saturation $\tau\left(t\right)=n\left(t\right)/N$.
We also remark that the Eq. (\ref{eq:tau-1}) holds for any $\alpha-$H�lder
urn function at least, and can be applied out of the box to more advanced
IR models that can still be embedded in the HLS urn model, like those
considered in \cite{Dosi Ermoliev,Dosi Kaniovsky}. For example, we
could have considered a market model where multiple products are present,
as far as we follow the market share of only one of them (say the
first one). If the customers follow the majority of the polled sample
with a probability $p$, and buy at random one of the $r>0$ available
products with probability $1-p$, the probability that the product
is purchased would have been
\begin{equation}
\pi_{k}^{*}\left(x\right):=p\,P_{k}\left(x\right)+\left(1-p\right)/r,
\end{equation}
although there are differences in the convergence properties (for
example this new function is not symmetric around $1/2$ and the autocorrelation
scaling of Ref. \cite{Kazuaki} may not hold) it is still possible
to repeat the same LD analysis, and find a similar phase structure
- this model will be discussed in detail elsewhere. Moreover, the
LDT techniques shown in Sections \ref{sec:Large-deviations} and \ref{sec:Zero-cost-trajectories}
goes beyond the HLS model, and could be adapted to find trajectories
of processes that are not directly embedded in the HLS model, such
as the one presented in the Ref. \cite{Dosi last}, where also the
possibility of losing customers is considered. It should be also possible
to extend the LDT to time-dependent urn functions that varies on a
time scale $O\left(N\right)$, maybe by considering a partition of
the range of $\tau$ into small subintervals where the urn function
can be approximated as a constant and then apply the equations given
before. We expect that even some quenched disordered versions of these
models may be studied, either by combining with the Replica Symmetry
Breaking theory \cite{PMV RSB} or the kernel methods of Ref. \cite{KERNEL THEO}.

\section{Trajectories of the DEK model\label{sec:Main-results}}

We explicitly write trajectories in closed form for the cases $k=1$
and $k=3$, distinguishing between those trajectories for which $u\left(0\right)\in\left[0,1\right]$,
the only possible for a virgin market start, from those crossing the
boundary values at some positive saturation (ie, $u\left(\tau\right)\in\left\{ 0,1\right\} $
for some $\tau>0$), that can have zero cost only in the $\tau_{0}>0$
case. For the  DEK model with $k=1$ we can evaluate the integral
that defines $\Pi_{1}$: 
\begin{equation}
\int_{u}^{y}\frac{d\alpha}{\pi_{1}\left(\alpha\right)-\alpha}=\frac{1}{1-p}\int_{u}^{x}\frac{d\alpha}{1-2\alpha}=\frac{1}{2\left(1-p\right)}\log\left(\frac{u-x_{0}}{x-x_{0}}\right),\label{eq:fdfdfd}
\end{equation}
this formula can be easily inverted, then from Eq. (\ref{eq:tau-1-2-1})
we find the equation for the trajectories 
\begin{equation}
u\left(\tau\right)=x_{0}+\left(x-x_{0}\right)/\tau^{\,2\left(1-p\right)},
\end{equation}
for $x\neq x_{0}$ these trajectories diverge for any $p>0$ when
$\tau\rightarrow0$, while for $x=x_{0}$ a unique non-divergent trajectory
exists for all $p$, and is $u\left(\tau\right)=x_{0}$. Notice that
in the limit of perfect trust $p=1$, that is equivalent to the classic
Polya Urn Model, each $u\left(\tau\right)=x$ becomes a non divergent
zero-cost trajectory for any share value $x$. For both late and early
starts, we find that the share $x_{N,n}$ always follows a single
zero-cost trajectory, that is therefore optimal. For a virgin market,
this trajectory is $u\left(\tau\right)=x_{0}$ and is independent
from the initial share, for late start we find the general convergence
point
\begin{equation}
x\left(u_{0},\tau_{0}\right)=x_{0}+\left(u_{0}-x_{0}\right)\tau_{0}^{\,2\left(1-p\right)}.
\end{equation}
This implies that for any initial condition, either early or late,
the entropy density of the DEK model with $k=1$ is zero only at the
critical value $x_{0}=1/2$ (and strictly negative otherwise) for
$p<1$, while for $p=1$ the entropy density is $\phi\left(x\right)=0$
at any point $x\in\left[0,1\right]$, as is expected for the Polya
Model \cite{Mahmoud}. 

\begin{figure}
\includegraphics[scale=0.29]{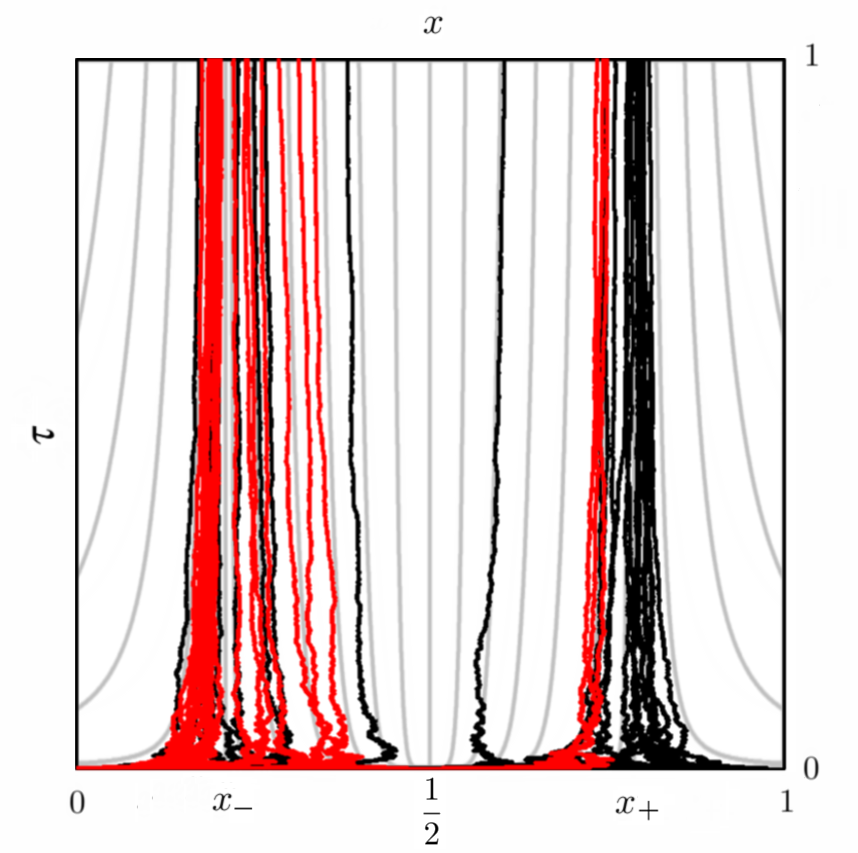}

\caption{\label{fig:DKE-model-for1}Zero-cost trajectories (the background
lines in light gray) of the DEK model for parameters $k=3$, $p=21/24>p_{c}$,
size $N=2^{14}$, early start at $\tau_{0}=2/N$, and initial condition
$u_{0}\in\left\{ 0,\,1/2\right\} $. For an early start in the market
at some negligible saturation $\tau_{0}=o\left(1\right)$ the trajectories
of the simulated process are scattered around the equilibrium points
$x_{+}$ and $x_{-}$ at the very beginning of the process and then
progressively stabilize on some zero-cost trajectory. The figure shows
$25$ realizations with initial conditions $X_{N,2}=0$ red (gray)
lines, $X_{N,2}=1$ black lines: notice that some processes starting
from zero where able to reach the nearby of $x_{+}$ anyway. }
\end{figure}

\begin{figure}
\includegraphics[scale=0.29]{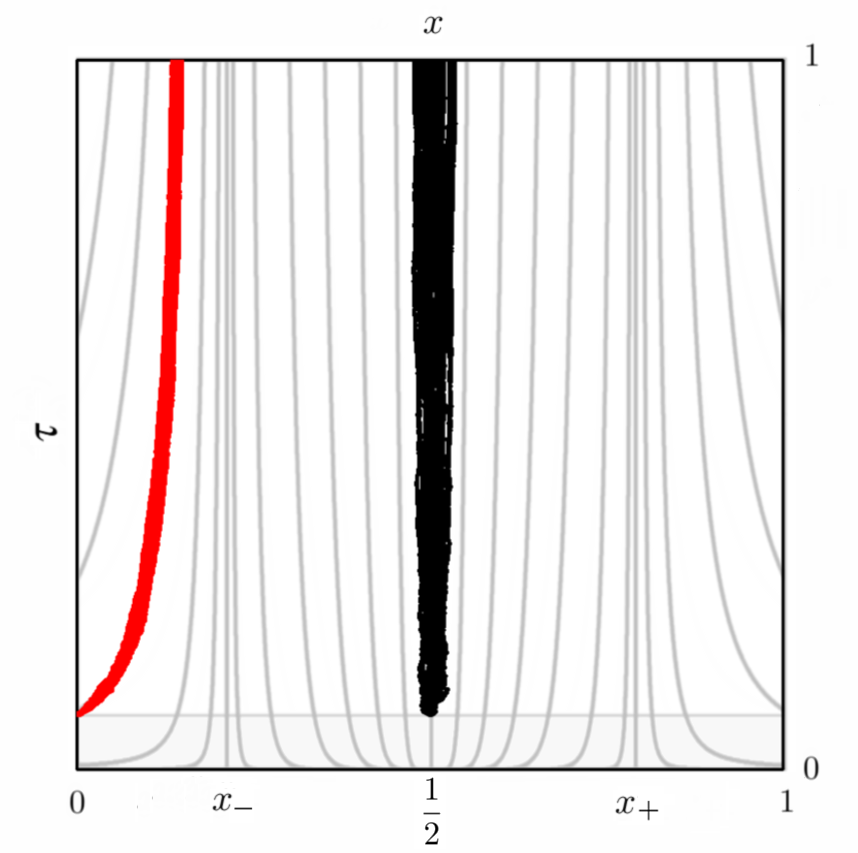}

\caption{\label{fig:DKE-model-for2}Zero-cost trajectories (background lines
in light gray, same of Figure \ref{fig:DKE-model-for1}) of the DEK
model for $k=3$, $p=21/24>p_{c}$, size $N=2^{14}$, late start at
$\tau_{0}=1/13$, and initial condition $u_{0}\in\left\{ 0,\,1/2\right\} $.
As one can see, in case of a late start at $\tau_{0}=O\left(1\right)$
the initial conditions become relevant, and the process follows a
single trajectory $u\left(\tau\right)$ emanating from the initial
condition at $(u_{0},\tau_{0})$. The figure shows $25$ realizations
with initial conditions $X_{N,2}=0$ red (gray) lines, $X_{N,2}=1$
black lines: the late start removes the degeneracy of the trajectories.}
\end{figure}

Most interesting to the IRT is the case $k=3$ with $p>p_{c}$, where
the lock-in phenomenon is possible: the general picture below $p_{c}=5/6$
is qualitatively the same that is found in the $k=1$ case, but above
$p_{c}$ and for virgin market start we observe a whole region $[x_{-},x_{+}]$
where the trajectories have sub linear entropy cost, although only
$u\left(\tau\right)=x_{\pm}$ are really optimal. Let compute the
trajectories for $k=3$: also in this case the integral can be evaluated
exactly, define the parameters 
\begin{equation}
\Delta:=6p-5,\ \ \Lambda^{2}:=\frac{\Delta}{4\left(2p-1\right)},\label{eq:PARA}
\end{equation}
the integral for $\Pi_{3}$ can be found via computer algebra,
\begin{equation}
\int_{u}^{x}\frac{d\alpha}{\pi_{3}\left(\alpha\right)-\alpha}=\frac{1}{\Delta}\log\left(\frac{1-\Lambda^{2}/\left(u-x_{0}\right)^{2}}{1-\Lambda^{2}/\left(x-x_{0}\right)^{2}}\right),\label{eq:paramagno}
\end{equation}
let introduce the $x-$dependent coefficient
\begin{equation}
\rho\left(x\right):=\Lambda^{2}/\left(x-x_{0}\right)^{2}-1,
\end{equation}
we can invert the Eq. (\ref{eq:paramagno}), and compute the equation
for the zero-cost trajectories also in the case $k=3$
\begin{equation}
u\left(\tau\right)=x_{0}\pm\sqrt{\frac{\Lambda^{2}}{1+\rho\left(x\right)/\tau^{\,\Delta}}},\label{eq:tau-2-1}
\end{equation}
where the plus and minus depends on weather the parameter $x$ lies
above or below $x_{0}=1/2$. Also in this case, for $x=x_{0}$ there
is a zero-cost trajectory for any $p$, in fact, for this value $\rho$
diverges, and the trajectory is $u\left(\tau\right)=1/2$. For $x\neq x_{0}$
we have to look weather the sign of $\Lambda^{2}$ is positive or
not, and we can see from Eq. (\ref{eq:PARA}) that when $p>p_{c}=5/6$
both $\Delta$ and $\Lambda^{2}$ are indeed positive quantities.
This implies that $\tau^{\Delta}$ converges to zero when also $\tau$
does, then the $u\left(\tau\right)-x_{0}$ converges to zero at the
admissible point $u\left(0\right)=1/2$. Finally, notice that $\rho\left(x\right)$
also must be positive, otherwise the formula inside the radical would
become negative for some $\tau>0$: then, any admissible trajectory
requires the further condition $x-x_{0}\in\left[-\Lambda,\,\Lambda\right],$
by confronting Eq. (\ref{eq:PARA}) with Eq. (\ref{eq:PARAK3}) we
can readily see that, as expected, $\Lambda$ is equal to half distance
between the convergence points $\left|x_{\pm}-x_{0}\right|$. Then,
the condition reduces to $x\in\left[x_{-},x_{+}\right]$, implying
that a non divergent zero-cost trajectory exists for any $x$ lying
between the convergence points. On the other hand, in the case of
a late start with an initial share $u\left(\tau_{0}\right)=u_{0}$
at some initial saturation $\tau_{0}>0$ there is always a unique
zero-cost trajectory emanating from $u_{0}$ and ending in 
\begin{equation}
x\left(u_{0},\tau_{0}\right)=x_{0}\pm\sqrt{\frac{\Lambda^{2}}{1+\rho\left(u_{0}\right)\tau_{0}^{\,\Delta}}},
\end{equation}
that is also optimal. In Figures (\ref{fig:DKE-model-for1}) and (\ref{fig:DKE-model-for2})
a simulation of the DEK model in the lock-in phase is shown for different
initial conditions, and confronted with its predicted trajectory.
See also Figure 1 and 2a of G. Dosi et al. (2019) \cite{Dosi last}
with the Figure 2.1 of Franchini (2017) \cite{Franchini} for the
zero-cost trajectories of the seminal model with $p=1$ by Arthur
et al. with virgin market initial conditions.

\part{Methods}

\section{Large deviations\label{sec:Large-deviations}}

The variational problem shown in Eq. (\ref{eq:THEVARIATION}) is deduced
from two central results of LDT, the Varadhan Integral Lemma and the
Mogulskii theorem (see the recent paper by Touchette \cite{Touchette}
for an introductory presentation, \cite{Pham} for some applications
to economy, or the very detailed book by A. Dembo and O. Zeitouni
\cite{Dembo Zeitouni} for a mathematical review). Now, instead of
considering the event $x_{N,N}=\left\lfloor xN\right\rfloor /N$,
let first study the simpler situation 
\begin{equation}
\Omega:=\{X_{N}\in\left\{ 0,1\right\} ^{N}:\,x_{N,N}\in\left[\alpha,\beta\right]\},
\end{equation}
where the sample paths end in the interval $\left[\alpha,\beta\right]$
that contains $x$. The limit entropy density of such event is 
\begin{equation}
\phi\left(\alpha,\beta\right):=-\lim_{N\rightarrow\infty}\frac{1}{N}\log P\left(X\in\Omega\right).
\end{equation}
The starting point is the formula for the probability mass of a sample
trajectory. Let 
\begin{equation}
Y_{N}=\{Y_{N,1},\,Y_{N,2},\,...\,,\,Y_{N,N}\}
\end{equation}
with $Y_{N,n}\in\{0,1\}$ be a possible path, hereafter \textit{sample
path}, of the process $X_{N}$, then, its probability mass $P\left(X_{N}=Y_{N}\right)$
according to the measure $P$ is given by the formula
\begin{equation}
P\left(X_{N}=Y_{N}\right)=\prod_{n\leq N}\pi\left({\textstyle y_{N,n}}\right)^{Y_{N,n}}\left(1-\pi\left(y_{N,n}\right)\right)^{1-Y_{N,n}}.
\end{equation}
From here we define the entropy density of the path:
\begin{equation}
S^{*}\left(Y_{N}\right):=-\frac{1}{N}\log P\left(X_{N}=Y_{N}\right),
\end{equation}
introducing the auxiliary function
\begin{equation}
H\left(\alpha,\beta\right):=\alpha\log\beta+\left(1-\alpha\right)\log\left(1-\beta\right)
\end{equation}
the entropy density before can be rewritten as
\begin{equation}
S^{*}\left(Y_{N}\right)=-\frac{1}{N}\sum_{n\leq N}H\left(Y_{N,n},\pi\left({\textstyle y_{N,n}}\right)\right).\label{eq:entropy}
\end{equation}
It will be useful to introduce a notation for the average respect
to the measure $P$
\begin{equation}
Ef\left(X_{N}\right):=\sum_{Y_{N}\in\left\{ 0,1\right\} ^{N}}P\left(X_{N}=Y_{N}\right)f\left(Y_{N}\right),
\end{equation}
in this notation the probability mass of the event $\Omega$ is
\begin{equation}
P\left(X_{N}\in\Omega\right)=EI\left(X_{N}\in\Omega\right)
\end{equation}
We perform a change of measure
\begin{equation}
P\left(X_{N}\in\Omega\right)=\sum_{Y_{N}\in\Omega}P\left(X_{N}=Y_{N}\right)=\sum_{Y_{N}\in\left\{ 0,1\right\} ^{N}}P\left(X_{N}=Y_{N}\right)I\left(Y_{N}\in\Omega\right)=\sum_{Y_{N}\in\left\{ 0,1\right\} ^{N}}e^{-NS^{*}\left(Y_{N}\right)}I\left(Y_{N}\in\Omega\right).\label{eq:fgfgfggf}
\end{equation}
such that the probability of $\Omega$ can be represented as follows:
\begin{equation}
P\left(X_{N}\in\Omega\right)=2^{N}E_{0}\,e^{-NS^{*}\left(Y_{N}\right)}I\left(Y_{N}\in\Omega\right),
\end{equation}
where $E_{0}$ is the average according to the uniform measure 
\begin{equation}
E_{0}f\left(X_{N}\right):=\frac{1}{2^{N}}\sum_{X_{N}\in\left\{ 0,1\right\} ^{N}}f\left(X_{N}\right)
\end{equation}
ie, the measure of a binary random walk.

The next step is to construct a continuous interpolation for the path
$Y_{N}$, we introduce the function 
\begin{equation}
\varphi:=\{\,\left(\left\lfloor \tau N\right\rfloor /N\right)\,y_{N,\left\lfloor \tau N\right\rfloor }+(\tau-\left\lfloor \tau N\right\rfloor /N)\,Y_{N,\left\lfloor \tau N\right\rfloor }:\,\tau\in\left[0,1\right]\},\label{eq:dfdfd-1-2}
\end{equation}
so that the probability of the sample path can be represented in terms
of $\varphi$. The interpolated trajectories are supported by 
\begin{equation}
Q\left(\Omega\right):=\{\,\varphi\in Q:\,Y_{N}\in\Omega\}.\label{eq:dfdfd-1-1}
\end{equation}

It can be shown that $S^{*}$ admits a continuous representation.
This representation can be informally derived by changing the sum
in Eq. (\ref{eq:entropy}) into an integral 
\begin{equation}
\frac{1}{N}\sum_{n\leq N}\rightarrow\int_{0}^{1}d\tau
\end{equation}
and apply the proper scaling to the arguments of $H$, i.e. 
\begin{equation}
Y_{N,n}\rightarrow\partial_{\tau}\varphi\left(\tau\right),\ \ \pi\left({\textstyle y_{N,n}}\right)\rightarrow\pi\left(\varphi\left(\tau\right)/\tau\right).
\end{equation}
Applying these substitutions we obtain the following entropy functional
that approximate $S^{*}$: 
\begin{equation}
S\left(\varphi\right):=-\int_{0}^{1}d\tau\,H\left(\partial_{\tau}\varphi\left(\tau\right),\pi\left(\varphi\left(\tau\right)/\tau\right)\right).
\end{equation}
It can be shown that if $\pi\in\left(0,1\right)$ this functional
is continuous respect to the sup norm
\begin{equation}
\left\Vert \varphi-\eta\right\Vert :=\sup_{\tau\in\left[0,1\right]}\left|\varphi\left(\tau\right)-\eta\left(\tau\right)\right|
\end{equation}
ie, is such that if $\varphi$ converges to $\eta$ in sup norm then
also 
\begin{equation}
\left|S\left(\varphi\right)-S\left(\eta\right)\right|\rightarrow0.
\end{equation}
In Ref. \cite{Franchini} it is actually shown that
\begin{equation}
\lim_{N\rightarrow\infty}\left|S^{*}\left(Y_{N}\right)-S\left(\varphi_{N}\right)\right|=0,
\end{equation}
then, if $S$ is continuous in the large $N$ limit holds
\begin{equation}
\log2-\phi\left(\alpha,\beta\right)=\lim_{N\rightarrow\infty}\frac{1}{N}\log E_{0}\,e^{-NS^{*}\left(Y_{N}\right)}I\left(Y_{N}\in\Omega\right)=\lim_{N\rightarrow\infty}\frac{1}{N}\log E_{0}\,e^{-NS\left(\varphi\right)}I\left(\varphi\in Q\left(\Omega\right)\right).\label{eq:fffff}
\end{equation}

This is enough to compute the rate function from Varadhan Integral
Lemma \cite{Dembo Zeitouni}. Informally, this theorem can be seen
as a rigorous functional version of the well known saddle-point method.
From Lemmas 4.3.2 and 4.3.4 of the book by Dembo and Zeitouni \cite{Dembo Zeitouni}
we obtain 
\begin{equation}
\lim_{N\rightarrow\infty}\frac{1}{N}\log E_{0}e^{-NS\left(\varphi\right)}I\left(\varphi\in Q\left(\Omega\right)\right)=-\inf_{\varphi\in Q\left(\alpha,\beta\right)}\left\{ S\left(\varphi\right)-S_{0}\left(\varphi\right)\right\} \label{eq:dfd}
\end{equation}
where $Q\left(\alpha,\beta\right)$ is the limit of the set $Q\left(\Omega\right)$,
i.e. 
\begin{equation}
\lim_{N\rightarrow\infty}Q\left(\Omega\right)=Q\left(\alpha,\beta\right)=\bigcup_{\gamma\in\left[\alpha,\beta\right]}Q\left(\gamma\right),
\end{equation}
and $S_{0}$ the rate function of a simple random walk with binary
steps, in our context would be the case $p=1/2$.

The rate function $S_{0}$ is provided by the Mogulskii Theorem \cite{Dembo Zeitouni}:
it states that the rate function of any process where the increments
form an i.i.d. sequence is given by 
\begin{equation}
S_{0}\left(\varphi\right)=-\int_{0}^{1}d\tau\,M\left(\partial_{\tau}\varphi\left(\tau\right)\right).\label{ratefunction-1-1}
\end{equation}
where $M$ is the Legendre transform 
\begin{equation}
M\left(\alpha\right):=\inf_{\beta\in\left[0,\infty\right)}\left\{ \alpha\beta-\zeta\left(\beta\right)\right\} ,
\end{equation}
of the moment generating function of the increments 
\begin{equation}
\zeta\left(\beta\right):=E_{0}\exp\left(\beta Y_{N,1}\right),
\end{equation}
in case of coin-flip distributed bi\-na\-ry variables:
\begin{equation}
E_{0}\exp\left(\beta Y_{N,1}\right)=\frac{1}{2}\sum_{Y_{N,1}\in\left\{ 0,1\right\} }\exp\left(\beta Y_{N,1}\right)=\frac{1+e^{\,\beta}}{2}.
\end{equation}
Applying the Legendre transform, and following the Mogulskii Theorem
\cite{Franchini,Dembo Zeitouni}, we find: 
\begin{equation}
S_{0}\left(\varphi\right)=-\log2+J\left(\varphi\right),
\end{equation}
where the functional $J$ is defined 
\begin{equation}
J\left(\varphi\right):=-\int_{0}^{1}d\tau H\left(\partial_{\tau}\varphi\left(\tau\right),\partial_{\tau}\varphi\left(\tau\right)\right),
\end{equation}
for any absolutely continuous $\varphi\in Q$, and is $-\infty$ otherwise,
i.e. those trajectories that are not absolutely continuous have zero
probability mass (and can be ignored). In the end it is found \cite{Franchini}
that the rate function is equal to
\begin{equation}
I\left(\varphi\right)=J\left(\varphi\right)-S\left(\varphi\right),
\end{equation}
Noticing that 
\begin{equation}
L\left(\alpha,\beta\right)=H\left(\alpha,\beta\right)-H\left(\alpha,\alpha\right)
\end{equation}
we arrive to the rate function as presented in Eq. (\ref{ratefunction}).

We remark that Eq. (\ref{eq:THEVARIATION}) cannot be deduced by contraction
principle, because the internal part of $Q\left(x\right)$ (the set
minus its boundary) is void, and then cannot be a continuity set for
the rate function $I\left(\varphi\right)$. Some additional arguments
would then be necessary to rigorously prove this result, where we
apply the contraction principle to the mass of $Q\left(\alpha,\beta\right)$,
and then show that is possible to take $\alpha,\beta\rightarrow x$. 

This proof is rather technical and we do not need to discuss it here,
the interested readers can find it in the proof section of the Ref.
\cite{Franchini}. Also, notice that the requirement that $\pi\in\left(0,1\right)$
is not fulfilled if $p=1$, in Ref. \cite{Franchini} a special surgery
on the set $Q$ is performed to a priori exclude the problematic trajectories
and extend the result to the general case $\pi\in\left[0,1\right]$.

\section{Scaling of the entropy inside the sub-linear region\label{sec:Scaling-of-the}}

\begin{figure}
\includegraphics[scale=0.29]{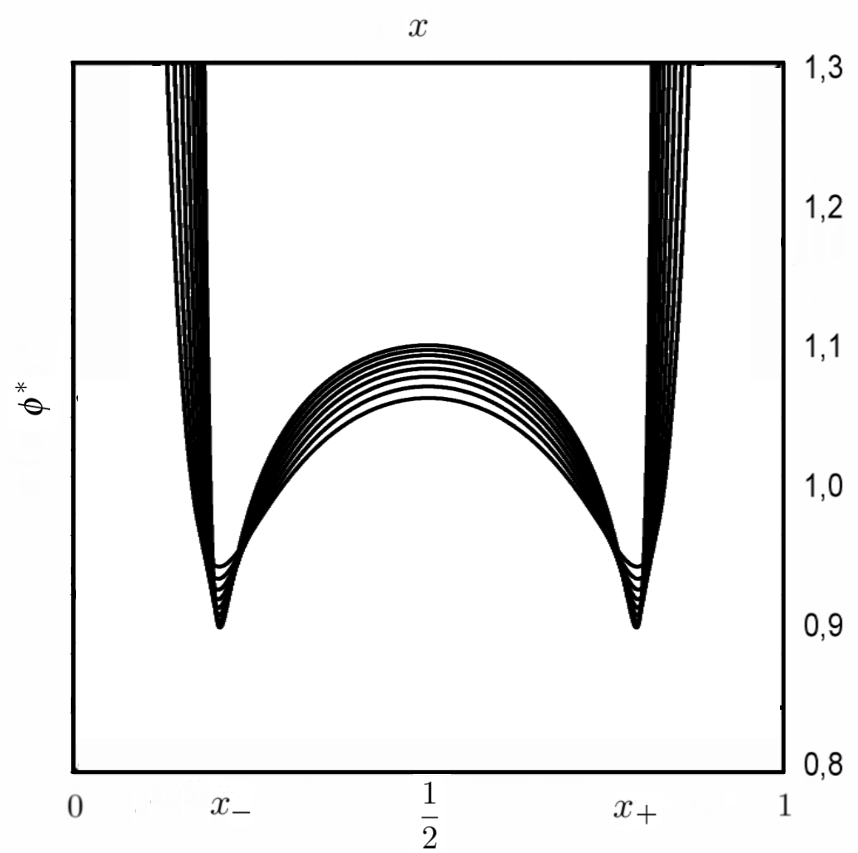}

\caption{\label{fig:Sub-linear-entropy}Sub-linear entropy scaling $\phi^{*}$
of the DEK model for $k=3$, $p=21/24>p_{c}$, size $N=2^{k}$ with
$10\leq k\leq17$, early start $\tau_{0}=2/N$, and initial condition
$u_{0}=1/2$ ($X_{N,2}=1$). The scaling limit is reached very slowly
in the sub-linear region, as the finite size corrections vanish only
logarithmically in the number of customers.}
\end{figure}

\begin{figure}
\includegraphics[scale=0.29]{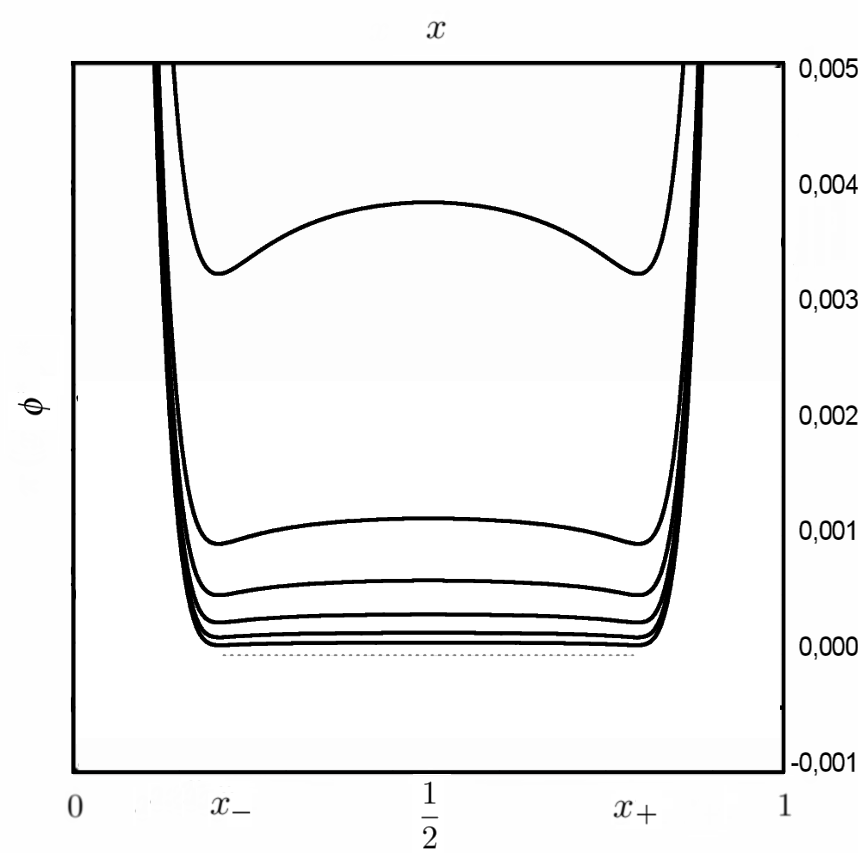}

\caption{\label{fig:Entropy-density-}Linear entropy scaling $\phi$ (entropy
density) of the DEK model for $k=3$, $p=21/24>p_{c}$, size $N=2^{k}$
with $12\leq k\leq17$, early start at $\tau_{0}=2/N$, and initial
condition $u_{0}=1/2$ ($X_{N,2}=1$). The scaling limit is reached
much faster in the linear regions. The gray dot line is the prediction
$\phi=0$ in the region between $\left[x_{-},\,x_{+}\right]$.}
\end{figure}

For a late market start at saturation $\tau_{0}>0$ and initial share
$u_{0}$ we have shown that the process follows a well defined trajectory
if $N$ is large enough, with a single convergence point $x\left(u_{0},\tau_{0}\right)$
where the $\phi$ is zero. The scaling of the entropy can be deduced
by noticing that this is the unique concentration point of the process,
then the probability mass of its nearby should be $O\left(1\right)$:
since for finite $N$ and any finite nearby of $x\left(u_{0},\tau_{0}\right)$
the mass must be distributed between a number of possible share values
that is of order $O\left(N\right)$, we expect that at the concentration
points the mass decays with some power of $N$. This reasoning suggests
that in the sub-linear region the entropy of the trajectory is logarithmic
in the number of potential customers, let define the sub-linear scaling
\begin{equation}
\phi^{*}\left(x\right):=-\lim_{N\rightarrow\infty}\,\frac{1}{\log N}\log P\left(x_{N,N}=\left\lfloor xN\right\rfloor /N\right)
\end{equation}
since there is a is unique concentration point, for any $\tau_{0}>0$
we can expect a monovariate probability mass function, then the predicted
sub-linear scaling would be divergent for any $x$ different from
$x\left(u_{0},\tau_{0}\right)$ and equal to some positive constant
otherwise. On the contrary, in the case of a virgin market start at
$\tau_{0}=0$ and for $p>p_{c}$ the limit of the entropy density
is found to be zero for any $x\in[x_{-},x_{+}]$, and therefore in
the lock-in phase the entropy of any trajectory that ends between
the points $x_{\pm}$ has a cost that is sub-linear in the potential
number of customers. In fact, is also possible to show \cite{Franchini}
that, for any continuous and invertible urn function, the limit $\phi\left(x\right)$
exists, it is strictly convex and negative from $x=0$ up to the first
point where the urn function crosses the diagonal, is zero from that
point to the last crossing, and then is convex negative again. A numeric
example of the scaling of $\phi$ and $\phi^{*}$ is in Figures (\ref{fig:Sub-linear-entropy})
and (\ref{fig:Entropy-density-}). 

Although the analysis of the zero-cost trajectories allows to establish
that in the early entry case the entropy in the region between the
convergence points is sub-linear, the exact scaling of the entropy
is not captured by this analysis. By the way, we remark that any deviation
from these trajectories on time scale $O(N)$ implies exponential
cost. Moreover, from the Corollary 6 of Ref. \cite{Franchini} follows
also the uniqueness of the solution for each $x\in\left(x_{1},x_{2}\right)$.
The probability mass current can flow along these trajectories only,
therefore, the current flowing through $(\varphi_{1},\varphi_{2})$
is a constant in $\tau$, 
\begin{equation}
P\left(\varphi\left(1\right)\in\left(x_{1},x_{2}\right)\right)=P\left(\varphi\left(\tau\right)\in\left(\varphi_{1}\left(\tau\right),\varphi_{2}\left(\tau\right)\right)\right),
\end{equation}
since can be also shown \cite{Franchini} that zero-cost trajectories
always emanate from the closest unstable equilibrium point, follows
that the entropy of the event $x_{N,N}\in\left(x_{1},x_{2}\right)$
should scale like the entropy near that point, that in this case is
$x_{0}=1/2$. It would be very interesting to have a general mathematical
theory that allows to find the exact rate at which the point $x_{0}$
expels its probability mass. An informal but general argument can
be found in Section III.B.2 of Jack (2019) \cite{Jack LD}. 

Interestingly, also Nakayama and Mori (2021) find that for urn functions
that are symmetric around $x_{0}$ the autocorrelation function satisfy
a universal logarithmic scaling for a suitable definition of the correlation
length. See Ref. \cite{Kazuaki} for further details.

\section{Cumulant generating function\label{sec:Cumulant-Generating-Function}}

We still didn't found much about the region $\phi\left(x\right)<0$:
it would be very interesting to have a method to compute the optimal
trajectories also in this region, perhaps this could be achieved by
properly deforming the zero-cost trajectories, or applying techniques
from Lagrange mechanics, or other optimal control methods. Although
this has not yet been achieved, we can still compute the shape of
$\phi$ outside the sub linear-region by analyzing the cumulant generating
function (CGF)
\begin{equation}
\xi\left(\lambda\right):=\lim_{N\rightarrow\infty}\frac{1}{N}\log\sum_{\Gamma\leq N}e^{-\lambda\Gamma}P\left(x_{N,N}=\Gamma/N\right),
\end{equation}
the right (left) behavior of $\phi\left(x\right)$ near the convergence
points can be deduced from the left (right) limit $\lambda\rightarrow0^{\pm}$
of the CGF before. Since the convergence points are always symmetric
around $x_{0}=1/2$ we only compute the limit from right. In Ref.
\cite{Franchini} is shown that, in general, the CGF satisfies the
following nonlinear differential equation at any $p$ and $k$ (eventually
any invertible $\pi$) 
\begin{equation}
\partial_{\lambda}\xi\left(\lambda\right)=\pi^{-1}\left(\frac{e^{\,\xi\left(\lambda\right)}-1}{e^{\,\lambda}-1}\right)\label{eq:equattons}
\end{equation}
with $\pi^{-1}$ inverse urn function, and we can study the behavior
at small lambda with a suitable perturbations theory (see next section).
The shape of $\phi$ near the convergence points is then computed
via the Legendre transform 
\begin{equation}
\phi\left(x\right)=\inf_{\lambda\in\left[0,\infty\right)}\left\{ \lambda x-\xi\left(\lambda\right)\right\} .\label{eq:laplkas}
\end{equation}
A possible informal derivation is as follows: let consider the difference
between the partition functions of the system at $N+1$ and $N$ customers
\begin{equation}
E\exp\,(\lambda Nx_{N,N}+\lambda X_{N+1,N+1})-E\exp\,(\lambda Nx_{N,N})={\textstyle \left(e^{\lambda}-1\right)}\,E(\pi\left(x_{N,N}\right)\exp\,(\lambda Nx_{N,N})),\label{eq:fdfd}
\end{equation}
consider the following equivalent expression for the CGF
\begin{equation}
\xi_{N}\left(\lambda\right)=\frac{1}{N}\log\,E\exp\,(\lambda Nx_{N,N}),
\end{equation}
define the auxiliary function
\begin{equation}
\delta_{N}\left(\lambda\right):=\left(N+1\right)\left(\xi_{N+1}\left(\lambda\right)-\xi_{N}\left(\lambda\right)\right),
\end{equation}
and the notation $E_{\lambda}$ for the tilted average, 
\begin{equation}
E_{\lambda}f\left(X_{N}\right):=\frac{Ef\left(X_{N}\right)\exp\left(\lambda Nx_{N,N}\right)}{E\exp\left(\lambda Nx_{N,N}\right)}.
\end{equation}
Using this notation and after some manipulations we arrive to the
identity 
\begin{equation}
\delta_{N}\left(\lambda\right)+\xi_{N}\left(\lambda\right)=\log\left(1+{\textstyle \left(e^{\lambda}-1\right)}\,E_{\lambda}\pi\left(x_{N,N}\right)\right),
\end{equation}
now we take the limit $N\rightarrow\infty$: from the existence of
$\phi$ follows that of $\xi$, then the limit of $\xi_{N}$ exists
and is 
\begin{equation}
\lim_{N\rightarrow\infty}\xi_{N}\left(\lambda\right)=\xi\left(\lambda\right),
\end{equation}
and can be shown \cite{Franchini} that, if the urn function $\pi$
is invertible, which is our case, then also its derivative exists,
and converges to $\partial_{\lambda}\xi$ in the limit 
\begin{equation}
\lim_{N\rightarrow\infty}E_{\lambda}x_{N,N}=\lim_{N\rightarrow\infty}\partial_{\lambda}\xi_{N}\left(\lambda\right)=\partial_{\lambda}\xi\left(\lambda\right).
\end{equation}
It is also possible to prove \cite{Franchini} that $x_{N,N}$ weakly
concentrates on its convergence point under the tilted average $E_{\lambda}$,
notice that the tilted average of $x_{N,N}$ is 
\begin{equation}
E_{\lambda}x_{N,N}=\partial_{\lambda}\xi_{N,N}\left(\lambda\right),
\end{equation}
therefore by weak convergence 
\begin{equation}
\lim_{N\rightarrow\infty}E_{\lambda}\pi\left(x_{N,N}\right)=\lim_{N\rightarrow\infty}\pi\left(E_{\lambda}x_{N,N}\right)=\pi\left(\partial_{\lambda}\xi\left(\lambda\right)\right).
\end{equation}
Finally, with a slightly more technical argument (see the proof section
of Ref. \cite{Franchini}) one can show that $\delta_{N}$ converges
to zero 
\begin{equation}
\lim_{N\rightarrow\infty}\delta_{N}\left(\lambda\right)=0,
\end{equation}
putting together we find 
\begin{equation}
\xi\left(\lambda\right)=\log\left(1+{\textstyle \left(e^{\lambda}-1\right)}\,\pi\left(\partial_{\lambda}\xi\left(\lambda\right)\right)\right),
\end{equation}
that is equivalent to Eq. (\ref{eq:equattons}). 

The DEK with $k=1$ is fully equivalent to the ERW, whose LDT properties
have been studied by Jack and Harris in both $p$ regimes \cite{Jack Harris}:
the urn function is 
\begin{equation}
\pi_{1}\left(x\right)=a+bx,
\end{equation}
with coefficients equal to
\begin{equation}
a=1-p,\ \ b=2p-1.
\end{equation}
From Eq. (\ref{eq:equattons}), the CGF satisfies the differential
equation 
\begin{equation}
a+b\,\partial_{\lambda}\,\xi\left(\lambda\right)=\frac{e^{\,\xi\left(\lambda\right)}-1}{e^{\,\lambda}-1},\label{eq:equattons-1-1-2}
\end{equation}
this equation can be integrated exactly by applying a proper substitution,
and then the Laplace method (see Section 3.3.2 Ref. \cite{Franchini}):
adapting the results from Corollary 10 Ref. \cite{Franchini} (see
also Jack and Harris \cite{Jack Harris}) we find that the CGF is
\begin{equation}
1-e^{-\xi\left(\lambda\right)}=\frac{a}{b}\,e^{\left(a/b\right)\lambda}\left({\textstyle 1-e^{-\lambda}}\right)^{1/b}\int_{1-e^{-\lambda}}^{1}dt\,\left(1-t\right)^{\left(a/b\right)-1}t^{-1/b}\label{eq:sssd}
\end{equation}
for $p>1/2$ and $\lambda>0$. Interestingly for $b>0$ ($p>1/2$)
the function is never analytic at $\lambda=0$, expanding for small
$\lambda$ we find a non vanishing term, of order $\lambda^{1/b}\log\lambda$
when $1/b$ is an integer number and $\lambda^{1/b}$ when is a real
number: derivatives of order higher than $1/b$ are singular at $\lambda=0$.

\section{Scaling of the master equation}

Numerically, we can study the shape of $\phi$ by computing the master
equation,
\begin{equation}
P\left(x_{N+1,N+1}=\Gamma/N\right)=\pi\left(\Gamma/N-1/N\right)P\left(x_{N,N}=\Gamma/N-1/N\right)+\left(1-\pi\left(\Gamma/N\right)\right)P\left(x_{N,N}=\Gamma/N\right),\label{eq:fffd}
\end{equation}
that can be integrated iteratively starting from the distribution
of the initial condition $x_{N,M}$. Notice that in practical numerical
tasks is not convenient to consider exponential quantities, and then
in our numerical tests we will consider the entropy
\begin{equation}
\Phi\left(\Gamma,N\right):=-\log P\left(x_{N,N}=\Gamma/N\right),
\end{equation}
in this form the master equation can be rewritten as follows
\begin{equation}
\exp\left(\Phi\left(\Gamma,N\right)-\Phi\left(\Gamma,N+1\right)\right)=\pi\left(\Gamma/N-1/N\right)\exp\left(\Phi\left(\Gamma,N\right)-\Phi\left(\Gamma-1,N\right)\right)+\left(1-\pi\left(\Gamma/N\right)\right).\label{eq:fffd-2}
\end{equation}

The Eq.s (\ref{eq:equattons}) and (\ref{eq:laplkas}) can be (informally)
deduced also from the master equation: in fact, the existence of $\phi$
suggests to try the following scaling 
\begin{equation}
\Phi\left(\Gamma,N\right)\rightarrow N\phi\left(\Gamma/N\right),
\end{equation}
that holds for large $N$. The left term of the master equation is
\begin{equation}
\Phi\left(\Gamma,N\right)-\Phi\left(\Gamma,N+1\right)\rightarrow-\left(N+1\right)\phi\left(\Gamma/\left(N+1\right)\right)+N\phi\left(\Gamma/N\right)\label{eq:dddd}
\end{equation}
while the right term is
\begin{equation}
\Phi\left(\Gamma,N\right)-\Phi\left(\Gamma-1,N\right)\rightarrow-N\phi\left(\Gamma/N\right)+N\phi\left(\left(\Gamma+1\right)/N\right),
\end{equation}
putting back into the master equation we find
\begin{multline}
\exp\left(\left(N+1\right)\phi\left(\Gamma/\left(N+1\right)\right)-N\phi\left(\Gamma/N\right)\right)=\\
=\pi\left(\Gamma/N-1/N\right)\exp\left(N\phi\left(\Gamma/N-1/N\right)-N\phi\left(\Gamma/N\right)\right)+\left(1-\pi\left(\Gamma/N\right)\right).\label{eq:fffd-2-1}
\end{multline}
Now, let apply the scaling $\Gamma/N\rightarrow x$ and $1/N\rightarrow dx$,
from this conditions we deduce that
\begin{equation}
\left(\Gamma+1\right)/N\rightarrow x+dx,
\end{equation}
\begin{equation}
\Gamma/\left(N+1\right)=\Gamma/N-\Gamma/\left(N\left(N-1\right)\right)\rightarrow x-xdx,
\end{equation}
the scaling of the entropy density is
\begin{equation}
\phi\left(\Gamma/\left(N+1\right)\right)\rightarrow\phi\left(x\right)
\end{equation}
\begin{equation}
\phi\left(\Gamma/\left(N+1\right)\right)\rightarrow\phi\left(x\right)-x\,\partial_{x}\phi\left(x\right)dx
\end{equation}
\begin{equation}
\phi\left(\left(\Gamma-1\right)/N\right)\rightarrow\phi\left(x\right)-\partial_{x}\phi\left(x\right)dx.
\end{equation}
In the end one obtains a non-linear differential equation
\begin{equation}
\pi\left(x\right)=\frac{\exp\,(x\,\partial_{x}\phi\left(x\right)-\phi\left(x\right))-1}{\exp\,(\partial_{x}\phi\left(x\right))-1},
\end{equation}
that reduces to Eq. (\ref{eq:equattons}) if one substitutes $\partial_{x}\phi\left(x\right)\rightarrow-\lambda$
and
\begin{equation}
x\,\partial_{x}\phi\left(x\right)-\phi\left(x\right)\rightarrow\xi\left(\lambda\right),\ \ \ x\rightarrow\partial_{\lambda}\xi\left(\lambda\right).
\end{equation}

It would be very interesting to have a general theory to solve these
differential equations for any $\pi$: at present, this can be done
only for linear urn functions.

\section{Perturbations theory\label{sec:Perturbations-theory} for $k=1$}

In these final sections we elaborate a first order perturbations theory
for the shape of $\phi$ outside the sublinear region. We find some
more critical values of the trust parameter $p$, that exist in both
the $k=1$ and $k=3$ cases. For $k=1$ only one $p^{*}$ exists,
beyond which the peak of the share distribution is not Gaussian anymore
(that is well known). Interestingly, in the case $k=3$ there are
two critical values: $p^{*}$, that is analogue to the case $k=1$,
and a $p^{**}$, beyond which the Gaussianity near the convergence
point seems restored, see Figures \ref{fig:3} and \ref{fig:4}.

To systematically understand the shape of $\phi$ it will be more
instructive to perform an approximate analysis. We put emphasis on
perturbation theory because is a simple method and does not require
special mathematical knowledge on ODE to be applied. We consider the
following general scaling at small $\lambda$ 
\begin{equation}
\xi\left(\lambda\right)\approx A\lambda+B\lambda^{2}+C\lambda^{\theta}
\end{equation}
where the approximate equality symbol $\approx$ is intended in the
sense that we are ignoring all terms of the kind $\lambda^{\theta}$
with $\theta>2$. This is because $\lambda$ is assumed to be small,
then the term $\lambda^{\theta}$ can rival with the regular terms
only if $\theta\leq2$, i.e., for $\theta>2$ the regular terms dominate
the first two moments of the distribution and $\lambda^{\theta}$
can be ignored. The derivative respect to $\lambda$ is 
\begin{equation}
\partial_{\lambda}\xi\left(\lambda\right)\approx A+2B\lambda+\theta C\lambda^{\theta-1}.
\end{equation}
Then, we approximate the right side of Eq. (\ref{eq:equattons-1-1-2}),
\begin{equation}
\frac{e^{\,\xi\left(\lambda\right)}-1}{e^{\,\lambda}-1}\approx A+\left(A\left(1-A\right)/2+B\right)\lambda+C\lambda^{\theta-1},
\end{equation}
equating the coefficients of the terms with equal power 
\begin{equation}
\left(\left(1-b\right)A-a\right)+\left(A\left(1-A\right)/2+\left(1-2b\right)B\right)\lambda+\left(1-b\theta\right)C\lambda^{\theta-1}\approx0,\label{eq:ddddf}
\end{equation}
we find the following values for $A$, $B$ and $\theta$:
\begin{equation}
A=\frac{a}{1-b}=\frac{1}{2},
\end{equation}
\begin{equation}
2B=-\frac{A\left(1-A\right)}{1-2b}=-\frac{1}{4\left(3-4p\right)},
\end{equation}
\begin{equation}
\theta=\frac{1}{b}=\frac{1}{2p-1},
\end{equation}
the amplitude $C$ is not captured by this expansion, and must be
determined in a different way, for example it could be obtained from
the exact expression of the CGF that is given before, but we don't
need it. 

We remark that, when $p>p^{*}=3/4$, i.e., when the derivative of
this urn function at the point of convergence $x_{0}$ goes above
$1/2$, then even the second order cumulant is super-linear, and the
shape $\phi\left(x\right)$ in the nearby of $x_{0}=1/2$ is not even
Gaussian anymore for $p\in\left(p^{*},1\right)$. This suggests some
phase change in the convergence mechanism of $x_{N,N}$: below $p^{*}$,
when the urn function derivative at the point $x_{0}$ is less than
$1/2$, we expect that $x_{N,N}$ will cross the critical value infinitely
many times in its evolution. But above the value $p_{c}$ the convergence
of $x_{N,N}$ has a slow down, according to an interesting mechanism
first described by Pemantle \cite{Pemantle Touch}, where $x_{N,N}$
approaches $x_{0}$ so slowly that it will never cross this point
(almost surely), and will accumulate in its neighborhood. 

The effects of this transition can be observed in the shape of $\phi$.
Let apply the Legendre transform to the expression of the CGF for
small $\lambda$, first we have to solve the equation 
\begin{equation}
x-\partial_{\lambda}\xi\left(\lambda\right)=0,
\end{equation}
inserting the approximation before we have
\begin{equation}
x-A-2B\lambda-\theta C\lambda^{\theta-1}\approx0.
\end{equation}
For $\theta\geq2$ the quadratic term is dominant at small $\lambda$,
and the previous condition reduces to
\begin{equation}
x-A-2B\lambda\approx0
\end{equation}
solving the equation we find the $\lambda$ that minimizes the Legendre
functional of Eq. (\ref{eq:laplkas})
\begin{equation}
\lambda\approx\frac{x-A}{2B},
\end{equation}
putting back in the expression for $\phi$ we find 
\begin{equation}
\phi\left(x\right)\approx B\left(\frac{x-A}{2B}\right)^{2}.
\end{equation}
If instead $1\leq\theta\leq2$ the quadratic term can be ignored in
favor of the non-linear term, the condition is 
\begin{equation}
x-A-\theta C\lambda^{\theta-1}\approx0
\end{equation}
The new condition bring to a different minimizer
\begin{equation}
\lambda\approx\left(\frac{x-A}{\theta C}\right)^{\frac{1}{\theta-1}},
\end{equation}
then we can compute the approximate shape, 
\begin{equation}
\phi\left(x\right)\approx\left(\theta C\right)^{-\frac{b}{1-b}}\left(\theta^{-1}-1\right)\left(x-A\right)^{\frac{1}{1-b}}.
\end{equation}
Summarizing, the shape of $\phi$ nearby the convergence point $y_{0}$
for the  DEK $k=1$ is approximately
\begin{equation}
\phi\left(x\right)\approx\begin{cases}
\begin{array}{l}
K_{0}\left|x-x_{0}\right|^{2}\\
K_{1}\left|x-x_{0}\right|^{1/\left(2-2p\right)}
\end{array} & \begin{array}{l}
0<p<p^{*}\\
p^{*}<p<1
\end{array}\end{cases},
\end{equation}
where the first constant is 
\begin{equation}
K_{0}=1/2B=2\left(3-4p\right)
\end{equation}
and $K_{1}$ must be determined from the exact form of $\xi$. To
keep the analysis simple we do not discuss the critical case $p=p_{c}$,
altough this also can be inferred from the exact form of $\xi$. 

\begin{figure}
\includegraphics[scale=0.29]{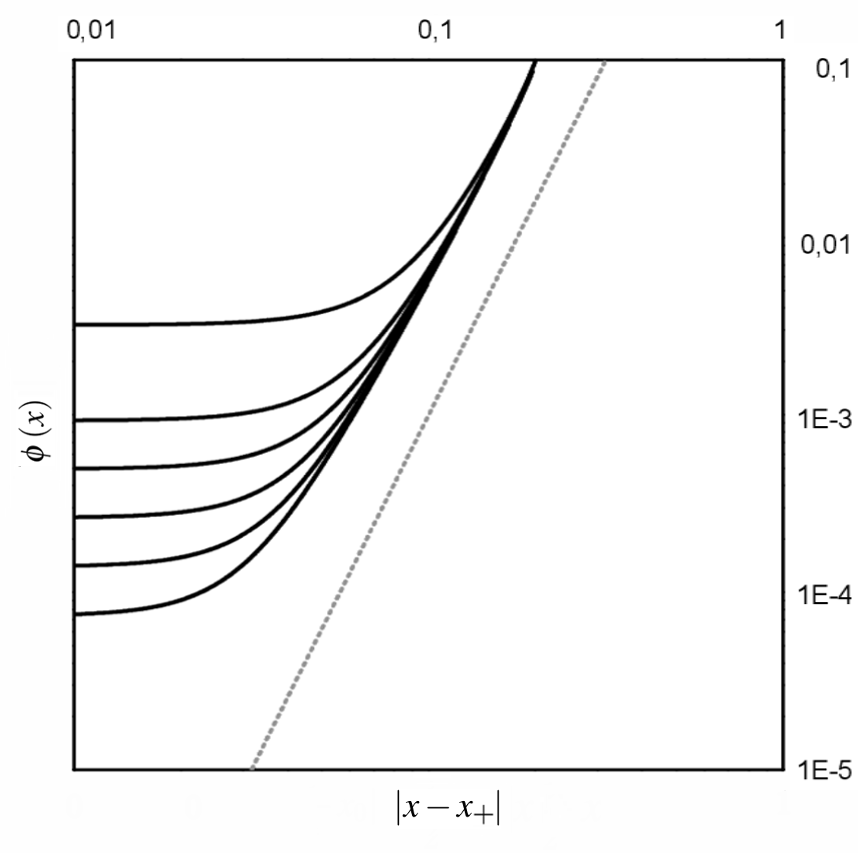}\caption{\label{fig:Test-of-perturbation1}Test of perturbation theory near
$x_{+}$ for $p_{c}<p<p^{**}$. The figure shows in double-logarithmic
scale the density $\phi$ of the DEK model near the convergence point
(in log-log scale) and the predicted behavior (gray dot line) $|x-x_{+}|^{4}$
from Eq. (\ref{eq:perturbationk3}) obtained via perturbations analysis.
The parameters are $k=3$, $p=21/24$, $N=2^{k}$ with $10\leq k\leq14$,
$\tau_{0}=2/N$ and $u_{0}=1/2$ ($X_{N,2}=1$). We remark that in
this range of $p$ the exact coefficient $K'_{2}$ cannot be determined
from perturbations, the gray dot line in the figure has been settled
to highlight the agreement of the predicted exponent with the simulation. }
\end{figure}

\section{Perturbations theory for $k=3$ \label{sec:Perturbations-theory-1}}

Concerning the generalized DEK $k=3$, its urn function is a third
degree polynomial of the kind
\begin{equation}
\pi_{3}\left(x\right)=a+cx^{2}-dx^{3}
\end{equation}
with null linear term, the other coefficients are 
\begin{equation}
a=1-p,\ \ c=3\left(2p-1\right),\ \ d=2\left(2p-1\right).
\end{equation}
The implicit differential equation for the CGF is
\begin{equation}
a+c\,\left(\partial_{\lambda}\xi\left(\lambda\right)\right)^{2}-d\,\left(\partial_{\lambda}\xi\left(\lambda\right)\right)^{3}=\frac{e^{\,\xi\left(\lambda\right)}-1}{e^{\,\lambda}-1},
\end{equation}
this equation cannot be solved (at best of our knowledge), but, by
looking at the behavior for small $\lambda$, we expect that below
$p_{c}$ the same picture of the linear case (with $k=1$) will arise,
although at different critical value $p^{*}=2/3$. This is because
the urn function has a flex at the convergence point $x_{0}$, i.e.
the urn function is locally linear. 

Let expand the urn function near the convergence point, for example
$x_{0}$, we can linearize it
\begin{equation}
\pi_{3}\left(x\right)\approx\pi_{3}\left(x_{0}\right)+\partial_{x}\pi_{3}\left(x_{0}\right)\left(x-x_{0}\right),
\end{equation}
where the derivative of $\pi_{3}$ is 
\begin{equation}
\partial_{x}\pi_{3}\left(x\right)=x\left(2c-3dx\right)=6\left(2p-1\right)x\left(1-x\right).
\end{equation}
We can use the results obtained for the linear case before, in the
region below $p_{c}$ we can take
\begin{equation}
b_{0}=\partial_{x}\pi_{3}\left(x_{0}\right)=6\left(2p-1\right)x_{0}\left(1-x_{0}\right)=\frac{3}{2}\left(2p-1\right)
\end{equation}
while above $p_{c}$we have
\begin{equation}
b_{\pm}=\partial_{x}\pi_{3}\left(x_{\pm}\right)=6\left(2p-1\right)x_{\pm}\left(1-x_{\pm}\right)=6\left(2p-1\right)\left(x_{0}+\Lambda\right)\left(x_{0}-\Lambda\right)=b_{0}\left(1-4\Lambda^{2}\right).\label{eq:ddds}
\end{equation}
Then, in the case $p<p_{c}$ we have $\Lambda=0$, recalling that
$\theta=1/b$ the sub-critical exponent is
\begin{equation}
\theta_{0}=\frac{2}{3\left(2p-1\right)},
\end{equation}
solving the equation $\theta=2$ we find $p^{*}=2/3$. For $p>p_{c}$
one has a positive $\Lambda$, substituting 
\begin{equation}
4\Lambda^{2}=\frac{6p-5}{2p-1}
\end{equation}
into Eq. (\ref{eq:ddds}) we find that 
\begin{equation}
\theta_{\pm}=\frac{1}{6\left(1-p\right)}
\end{equation}
therefore, there is another critical point where the shape of $\phi$
changes again, solving the equation we find it at 
\begin{equation}
p^{**}=11/12.
\end{equation}

Then, above $p_{c}$, another special trust parameter $p^{**}$ can
be identified, that corresponds to the value at which the derivative
of the urn function near $x_{\pm}$ (that above $p_{c}$ is decreasing
in $p$) goes once again below $1/2$. We predict that in this last
region the convergence mechanism below $p^{*}$ is restored, although
with a different convergence point. 

Putting these considerations together, we find the following approximate
shape of $\phi$ for $k=3$, 
\begin{equation}
\phi\left(x\right)\approx\begin{cases}
\begin{array}{l}
K'_{0}\left|x-x_{0}\right|^{2}\\
K'_{1}\left|x-x_{0}\right|^{\frac{2}{5-6p}}\\
K'_{2}\left|x-x_{\pm}\right|^{\frac{1}{6p-5}}\\
K'_{3}\left|x-x_{\pm}\right|^{2}
\end{array} & \begin{array}{l}
0<p<p^{*}\\
p^{*}<p<p_{c}\\
p_{c}<p<p^{**}\\
p^{**}<p<1
\end{array}\end{cases}\label{eq:perturbationk3}
\end{equation}
where we implicitly assumed that $\left|x\right|>1/2+\Lambda$, since
in the region between the convergence points we already know that
$\phi=0$. The constants $K'_{0}$ and $K'_{3}$ are computed like
in the $k=1$ case, one finds 
\begin{equation}
K'_{0}=2\left(1-2b_{0}\right)=4\left(3-2p\right),
\end{equation}
in the region below $p^{*}$ and 
\begin{equation}
K'_{3}=\frac{2\left(1-2b_{\pm}\right)}{1-4\Lambda^{2}}=\frac{1+12\left(1-2p\right)\left(1-p\right)}{2\left(1-p\right)}
\end{equation}
in the region above $p^{**}.$ On the contrary, the constants $K'_{1}$
and $K'_{2}$ cannot be determined using perturbations, and should
be found by other methods. A numerical check of the exponent in the
region $p_{c}<p<p^{**}$ is in Figure (\ref{fig:Test-of-perturbation1}).

\section*{Acknowledgments}

We thank Giovanni Dosi (Scuola Superiore Sant'Anna) and two anonymous
referees of Physical Review E for their useful comments. This research
has received funding from European Research Council (ERC) under the
European Union's Horizon 2020 research and innovation programme (Grant
Agreement No {[}694925{]}).

\end{document}